\documentclass[12pt]{article}
\usepackage{graphics,graphicx,amscd,amssymb,amsmath,color}
\begin{document}
\begin{center}
{\Large {A variant $\beta$-Wythoff Nim on Beatty's theorem}}
\bigskip
\newline Urban Larsson, 18 March 2011
\end{center}
\bigskip
\begin{center}{\bf Abstract }\end{center} \vspace{-0.3cm}
We give short rules for two-pile take-away games satisfying 
that a pair of complementary homogeneous Beatty sequences 
together with $(0,0)$ constitute a complete set of $P$-positions.
\bigskip
\bigskip

\noindent {\sc A new construction strictly in between Nim and $k$-Wythoff Nim.} 
Let us recall the rules of $k$-Wythoff Nim [F], $k$ a positive integer. 
The available positions are $(x,y)$, $x$ and $y$ non-negative integers. 
The legal moves are\\\\ 
\noindent (1A) Nim type: $(x,y)\rightarrow (x-t,y)$, if $x-t\ge 0$ and 
$(x,y)\rightarrow (x,y-t)$, if $y-t\ge 0$.\\\\
\noindent (1B) Extended diagonal type: 
$(x, y)\rightarrow (x - s,y - t)$ if $\mid s - t \mid < k$ and 
$x-s\ge 0, y-t\ge 0$.
\\\\
Hence, this game is a so-called impartial take-away game [WW]. 
We play normal play, that is the last player to move wins. By the 
rules (1A) and (1B) we note that  
$k$-Wythoff Nim is a so-called `invariant' [DR, LHF] take-away game, that is, 
each available move is legal from any position as long as the 
resulting position has non-negative coordinates.

In this note we study another type of take away games, where certain 
positions have some local restrictions on the set of otherwise 
`invariant' moves. Such games are sometimes called `variant', eg. [DR, LHF].
\\\\
\noindent {\bf Example 1} As usual, let $\pi = 3,14\ldots$ denote the ratio of 
the circumference of a circle to its diameter. Let the game rules 
of $\pi$-Wythoff Nim (our notation) be as in $k$-Wythoff Nim 
with $k = \lfloor \pi \rfloor = 3$, except that 
if a player plays from a position where one of the coordinates 
equals $\lfloor \pi n\rfloor$, for some $n\in \mathbb{N}$, then only Nim 
type moves (1A) are allowed. This latter rule clearly 
makes $\pi$-Wythoff Nim a `variant' game. Denote with 
$\alpha = \pi /(\pi -1)$. Then clearly both $\alpha$ and $\pi$ are irrational, 
so that, by Beatty's/Rayleigh's [B,R] theorem, the sequences 
$(\lfloor n\alpha \rfloor)$ and $(\lfloor n\pi \rfloor)$ 
are complementary on the positive integers. The main theorem of this note 
says that the $P$-positions of $\pi$-Wythoff Nim are identical to 
the set 
$\{(\lfloor n\alpha \rfloor, \lfloor n\pi \rfloor), (\lfloor n\pi \rfloor, 
\lfloor n\alpha \rfloor)\mid n\in \mathbb{Z}_{\ge 0}\}$.
\\\\

In general, fix an irrational  $2<\beta$. Then play the following variant 
game on the pairs of non-negative integers:
The moves are as in  $k$-Wythoff Nim with $k = \lfloor \beta \rfloor$ (1A) 
and (1B), except if one of the coordinates is of the form 
$\lfloor \beta n\rfloor$, 
$n\in \mathbb{Z}_{>0}$, then only Nim-type moves (1A) are allowed. Denote this 
game by $\beta$-Wythoff Nim.\\\\

\noindent{\bf Main Theorem} The $P$-positions of $\beta$-Wythoff Nim are 
$$\{(\lfloor n\alpha \rfloor, \lfloor n\beta \rfloor), (\lfloor n\beta \rfloor, 
\lfloor n\alpha \rfloor)\mid n\in \mathbb{Z}_{\ge 0}\},$$ where $n$ ranges 
over the non-negative integers and where $\alpha=\beta/(\beta-1)$.\\\\

\noindent{\bf Proof.} $P\rightarrow N$: We have to prove that from 
each position of 
the form 
\begin{align}\label{h}
(\lfloor \alpha n\rfloor, \lfloor\beta n\rfloor ), 
\end{align}
there is no move to a position of the same form, or to its 
symmetric counterpart of the form 
$(\lfloor \beta n\rfloor, \lfloor\alpha n\rfloor )$. 
So, suppose that we play from a position of the form in (1). Then, 
by the rules of game, only (1A) Nim type moves are allowed so that, 
by complementarity [B,R], there is no move to a position of the same form.\\\\ 
\noindent $N\rightarrow P$: If the candidate $N$-position 
$(x, y)$, $x\le y$, has a coordinate of the form $\lfloor \beta n\rfloor$ 
then we have to show that a Nim type (1A) move suffices. 
If $x = \lfloor \beta n\rfloor$ then 
move $(x,y)\rightarrow (\lfloor \beta n\rfloor, \lfloor \alpha n\rfloor)$.  
If $y = \lfloor \beta n\rfloor$ and $x > \lfloor \alpha n\rfloor$ move 
$(x,y)\rightarrow (\lfloor \alpha n\rfloor, \lfloor \beta n\rfloor)$.  
If $y = \lfloor \beta n\rfloor$ and $x < \lfloor \alpha n\rfloor$, by 
complementarity [B,R], there is an $m < n$ such that there is 
a Nim type (1A) move of precisely one of the forms  
$(x,y)\rightarrow (\lfloor \alpha m\rfloor, \lfloor \beta m\rfloor)$ or 
$(x,y)\rightarrow (\lfloor \beta m\rfloor, \lfloor \alpha m\rfloor)$. 
(We use that $\lfloor \alpha m\rfloor \le \lfloor \beta m\rfloor < \lfloor \beta n\rfloor$.) 

Otherwise, a Nim type move obviously suffices if 
$\lfloor \alpha m\rfloor = y > \lfloor \beta n\rfloor$ for some $n<m$ 
and $x = \lfloor \alpha n \rfloor$, so suppose that 
\begin{align}\label{g}
\lfloor \alpha m\rfloor = y < \lfloor \beta n\rfloor \text{ and }
x = \lfloor \alpha n \rfloor. 
\end{align}
(Still with $m>n$.) For this case, by complementarity, 
a Nim type move to a candidate $P$-position does not exist, 
so we have to find a (1B) ``extended diagonal'' type move.

Note that an ``ordered vector subtraction'' of consecutive 
candidate $P$-positions gives an expression of the form 
$$(\lfloor\alpha n\rfloor,\lfloor\beta n\rfloor) - 
(\lfloor\alpha (n-1)\rfloor,\lfloor\beta (n-1)\rfloor),$$ 
$n\in \mathbb{Z}_{>0}$, which equals precisely one of 
the four ordered pairs of differences: 
$$(1, \lfloor \beta \rfloor),$$
$$(1, \lfloor \beta +1\rfloor),$$
$$(2, \lfloor \beta \rfloor ),$$ or
$$(2,\lfloor \beta + 1\rfloor).$$
The difference of the coordinates in such a pair of 
differences is bounded by $\pm \lfloor \beta \rfloor$. By our assumption (2)  
this gives that there is a type (1B) move to a position of the form of a 
candidate $P$-position, 
$(\lfloor \alpha p\rfloor, \lfloor \beta p \rfloor)$ , $p<n$. 
In fact, a ``worst case scenario'' would be from an $N$-position of the form 
$(x,y) = (\lfloor\alpha n\rfloor+t,\lfloor\beta n\rfloor-1+t)$, 
$t\in \mathbb{Z}_{\ge 0}$, together with the above second case difference pair, 
$(1, \lfloor \beta +1\rfloor).$ But, indeed, here a move of type (1B) 
suffices to the $P$-position 
$(\lfloor\alpha (n-1)\rfloor, \lfloor\beta (n-1)\rfloor)$.
\hfill $\Box$
\bigskip
\bigskip

\noindent{\sc Questions and remarks.} Suppose that we fix a $\beta$ and then 
increase the density of the pairs of sequences
from 1 to say an arbitrary number $\gamma > 1$ (or decreases to 
a density $<1$) where $\alpha$ is defined via $1/\alpha + 1/\beta = \gamma$. 
Given candidate $P$-positions as above (Main Theorem), is there 
still a ``succinct'' and non-trivial  way of formulating the game rules without 
revealing both irrationals or/and the joint density of the sequences? 
As a remark, observe that neither $\alpha$ nor the 
density 1 is given away in the presentation of the rules of $\beta$-Wythoff 
Nim. In [LHF] invariant 
game rules are given for candidate $P$-positions constructed from 
complementary Beatty sequences, but not in a single case have we found 
a ``succinct'' description. In this note we have chosen to remove 
the nice condition of invariance from 
the game rules and, maybe even more notably, one of the coordinates of the 
candidate $P$-positions is revealed within the game rules. 
This could be argued to be a severe drawback in a definition of the rules 
of a game. But, on the other hand, we were able to give a very succinct 
formulation, without a complete trivialization of game rules, for all 
complementary Beatty sequences and these are uncountably many.
\\\\
{\bf Bibliography}

\noindent [B] S. Beatty, Problem 3173, \emph{Amer. Math. Monthly},
{\bf 33} (1926) 159.

\noindent [DR] E. Duch\^{e}ne and M. Rigo, Invariant Games,
\emph{Theoretical Computer Science}, {\bf 411} 34-36 (2010), 
pp. 3169-3180 

\noindent [F] A.S. Fraenkel, How to beat your Wythoff games' 
opponent on three fronts, \emph{Amer. Math. Monthly} {\bf 89} (1982) 353-361.

\noindent [LHF] U. Larsson, P. Hegarty, A. S. Fraenkel, 
Invariant and dual subtraction games resolving 
the Duch\^{e}ne-Rigo conjecture \emph{Theo. Comp. Sci.} {\bf 412} (2011) 729–735.

\noindent [R] J. W. Rayleigh. The Theory of Sound,\emph{Macmillan, London}, 
(1894) p. 122-123.

\noindent [WW] E. R. Berlekamp, J. H. Conway, R.K. Guy,
\emph{Winning ways}, {\bf 1-2} Academic
Press, London (1982). Second edition, {\bf 1-4}.
A. K. Peters, Wellesley/MA (2001/03/03/04).

\end{document}